# QUELQUES CONDITIONS POUR UN BON DÉPART EN RÉSOLUTION DE PROBLÈMES À L'ÉCOLE ÉLÉMENTAIRE *PRIMUM NON NOCERE*

Magali HERSANT[1]

INSPÉ Académie de Nantes, CREN, Nantes Université

Yves THOMAS[2]

Ancien formateur IUFM

**Résumé :** Les élèves de l'école élémentaire adoptent parfois face à un problème arithmétique des comportements que l'on déplore. Ces derniers peuvent résulter de certaines pratiques pédagogiques et pourraient être atténués par quelques actions assez faciles à mettre en place en classe. Nous proposons quatre pistes portant sur le choix et la formulation des problèmes et l'explicitation d'éléments du contrat didactique.

*Mots-clés : contrat didactique, résolution de problèmes, cycle 2, validation, calcul.*

## INTRODUCTION

L'apprentissage de la résolution de problèmes mathématiques à l'école élémentaire est source de nombreuses difficultés pour les élèves et de nombreux questionnements pour les enseignants. En particulier, dans le domaine de la résolution de problèmes arithmétiques textuels, les élèves peuvent : répondre au hasard, se baser sur des indices de surface non pertinents (mots de l'énoncé, opération utilisée dans le problème précédent ou le sujet de la leçon du jour), manquer d'esprit critique par rapport aux résultats obtenus (Verschaffel et al., 2000) et développer des règles d'actions erronées en lien avec le contrat didactique qu'ils imaginent (il faut utiliser tous les nombres de l'énoncé, il faut faire une soustraction parce que l'énoncé contient le mot « enlever » ou parce que la leçon de la veille portait sur la soustraction..).

Ces difficultés sont souvent mises en relation avec la difficulté des élèves à construire une représentation du problème. (Julo, 1995 ; Fagnant & Vlassis, 2013 ; Thevenot & Barouillet, 2015 ; Fagnant, 2019, 2022). L'objectif d'aider les élèves à construire une représentation des problèmes conduit actuellement les chercheurs," , mais aussi les décideurs, à formuler des propositions sur les différentes représentations, schématisations et modélisations des problèmes proposés aux jeunes élèves. Ces propositions traitent de la manière de résoudre un problème. Elles supposent implicitement que, dans la classe, l'activité de résolution de problème obéit plus ou moins à un contrat didactique que l'on peut formuler de la façon suivante :*résoudre un problème numérique, c'est apporter une réponse vraie à une question portant sur une situation où des nombres interviennent. La réponse est élaborée en utilisant les informations fournies, les connaissances que l'on a sur les nombres (faits numériques ou propriétés plus générales), et une réflexion rationnelle*. Ce contrat nous semble pertinent mais les pratiques « ordinaires » s'en écartent sur certains points. Souvent, en classe, l'attente est que la réponse soit « juste », ou « bonne » (conforme à ce qu'attend l'enseignant), plutôt que « vraie » (conforme à une réalité objective observable). Par ailleurs, des exigences extérieures à ce contrat sont souvent posées, par exemple « il faut faire des opérations ». De notre point de vue, ce malentendu fondamental explique les difficultés de certains élèves dans le domaine de la résolution de problèmes numériques et ce contrat doit donc être mis en place en amont du travail sur la façon de se représenter un problème pour le résoudre.

---

[1] magali.hersant@univ-nantes.fr

[2] primaths-yves-thomas@orange.fr

Dans cet article, nous proposons quatre conditions didactiques qui visent à installer un tel contrat au tout début de l'école élémentaire. Les deux premières sont étroitement liées aux choix matériels et à la formulation des problèmes. Il s'agit, par opposition à l'évaluation par l'enseignant (Margolinas, 1993), de permettre la validation par le matériel, en lien avec ce que les didacticiens nomment les rétroactions du milieu pour l'apprentissage (Brousseau, 1998). Il s'agit aussi de proposer des problèmes sans texte ou avec très peu de texte pour ne pas pénaliser les petits lecteurs, et ainsi contribuer à la dévolution des problèmes (Brousseau, 1998). Les deux autres conditions concernent plutôt les « règles » de production des réponses dans la classe et les attentes de l'enseignant. Il s'agit de permettre systématiquement de répondre « je ne sais pas » et d'indiquer dès les premiers problèmes qu'il faut répondre sans dénombrer.

Dans la première partie de l'article nous proposons douze problèmes pour le CP et le CE1 faciles à mettre en œuvre avec des objets présents dans la classe[3]. Nous montrons ensuite en quoi ils permettent l'installation du contrat que nous visons.

## I. DOUZE PROBLÈMES DU DÉBUT DU CP AU MILIEU DU CE1

Pour des raisons de présentation, les problèmes sont classés en fonction du matériel utilisé comme support, et non par ordre de difficulté ou par ordre chronologique d'utilisation en classe. Par ailleurs, la présence de matériel ne signifie pas que les élèves manipulent. En effet, les boites des problèmes 1 et 2 ne sont ouvertes qu'après la phase de recherche. De la même façon, la bande inconnue n'est retournée (problème 3) et les grandes bandes ne sont fabriquées (problème 4) qu'après la recherche. Cela impose le plus souvent que le matériel soit seulement collectif.

Ces douze problèmes ne suffisent évidemment pas pour une année et demie de travail, mais ils permettent d'illustrer quelques choix importants, qui seront discutés dans la suite de ce texte.

### 1. Avec une boite et des billes

L'enseignant décrit le contenu de la boite et fait au tableau les dessins ci-dessous pour que les données ne soient pas oubliées[4].

---

[3] La facilité de mise en œuvre n'est pas l'objet de l'article mais elle est importante pour la prise en compte de nos propositions par les enseignants. C'est pourquoi nous la mentionnons ici.

[4] Les dessins effectués au tableau respectent les couleurs annoncées dans l'énoncé. Par exemple pour le problème 1 le trois est écrit en rouge, le point d'interrogation en bleu et le neuf en noir.

Problème 1

Dans cette boite, j'ai mis 9 billes. Il y a 3 billes rouges, les autres sont bleues. Combien de billes bleues ?

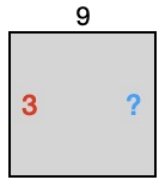

Problème 2

Dans cette boite, j'ai mis 5 billes rouges et 4 billes bleues. Je vais compter toutes les billes. Combien vais-je trouver de billes ?

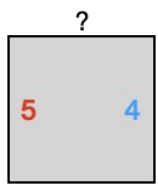

## 2. Avec des bandes quadrillées

Le matériel utilisé est constitué de bandes affichées au tableau. La face cachée est quadrillée. Le nombre écrit sur la face visible indique le nombre de cases au dos.

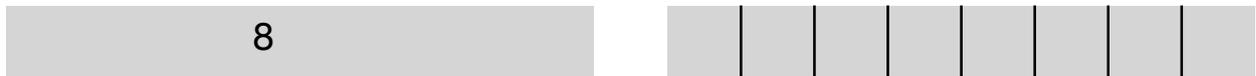

Problème 3

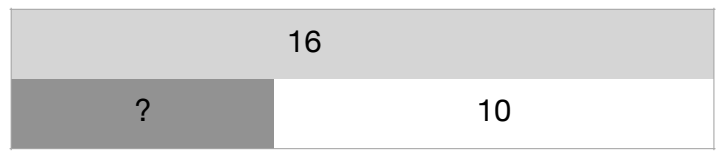

Problème 4

Je vais faire deux grandes bandes, une noire et une grise. Quelle sera la plus longue ?

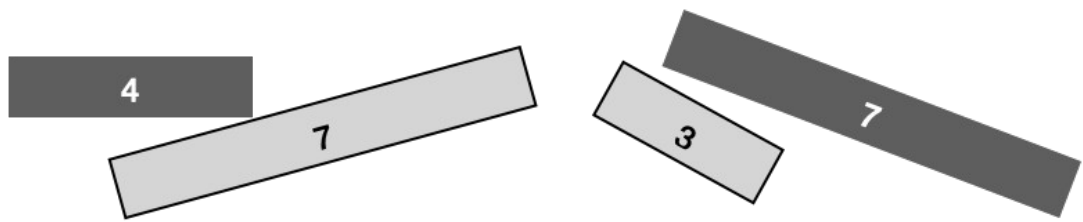

Il nous semble important de formuler dès à présent deux autres remarques au sujet de ces problèmes. D'abord, malgré une certaine similitude graphique, les problèmes de bandes que nous proposons n'ont de lien ni avec les réglettes Cuisenaire© ni avec les schémas en barre. L'identité de couleur des bandes de même longueur dans certains problèmes est locale et temporaire, il s'agit seulement d'insister sur le fait que, dans ce problème, les bandes inconnues sont identiques. Contrairement aux réglettes Cuisenaire©, la couleur n'est pas associée à un nombre. Les élèves ne peuvent pas se référer à des connaissances intimement liées au matériel, mais sans signification mathématique hors de ce matériel comme : « la réglette bleue est toujours aussi longue que la jaune et la violette mises bout à bout ». Ensuite, les bandes que nous utilisons ne sont pas des représentations d'un nombre de billes, de gâteaux ou d'euros comme c'est le cas quand on utilise des schémas en barre, ce sont des objets matériels à propos desquels on se questionne.

## 3. Avec des cartes à points

Le matériel est constitué de cartes affichées au tableau. Elles portent au verso des points disposés en privilégiant les constellations du dé. Le recto de chaque carte indique le nombre de points dessinés au verso.

Problème 5

Je vais retourner les cartes et compter tous les points. Combien y a-t-il de points en tout ?

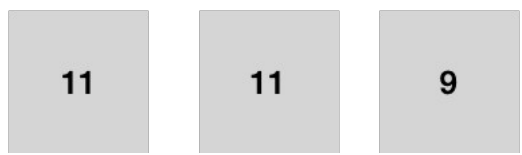

Problème 6

J'ai regardé derrière les cartes et j'ai compté tous les points. Il y a 17 points en tout.

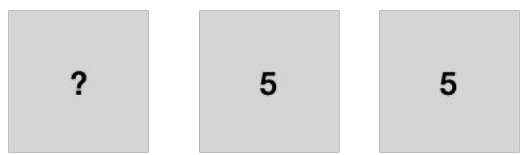

## 4. Avec une balance de Roberval et des objets de masse identique

Les objets peuvent être éventuellement de couleur différente, par exemple des pièces identiques d'un jeu de construction type « Duplo ».

Par convention, si l'on écrit 7 à côté d'un plateau du dessin de la balance, cela signifie qu'on placera une pile de 7 pièces sur ce plateau[5].

---

[5] Les problèmes 7, 8, 22 et 23 demandent d'effectuer des comparaisons, mais ce ne sont pas des problèmes de comparaison au sens de Vergnaud : il ne s'agit pas de déterminer un nombre.

Problème 7

De quel côté la balance va-t-elle pencher ?

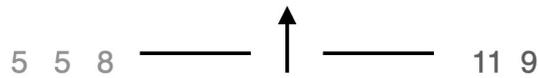

5 5 8          11 9

Problème 8

Je veux que la balance soit en équilibre en ajoutant une pile de pièces. Que faut-il ajouter ?

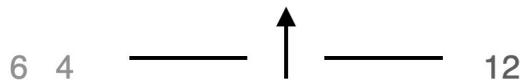

6 4          12

## 5. Avec de la monnaie

Le matériel est constitué de pièces fictives d'un euro et de quelques objets à vendre, dont le prix est affiché. L'enseignant ou un élève effectue les achats devant la classe.

Problème 9

J'ai 12 €. J'achète un gâteau à 5 € et un pain à 2 €. Maintenant, combien ai-je dans mon porte-monnaie ?

Problème 10

J'ai 40 €. Ai-je assez d'argent pour acheter un livre qui coûte 18 € et un CD qui coûte 15€ ?

## 6. Avec des rectangles quadrillés

Le matériel est constitué de rectangles affichés au tableau, dont la face cachée est quadrillée alors que la face visible est unie. Différentes indications peuvent être données : le nombre de carreaux sur un côté du rectangle ou le nombre total de carreaux à l'intérieur du rectangle.

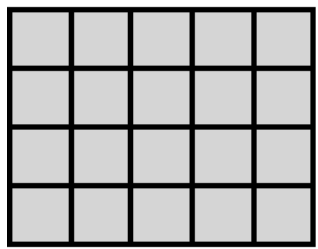 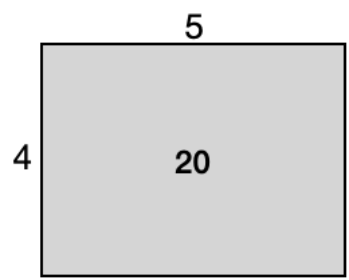

Problème 11

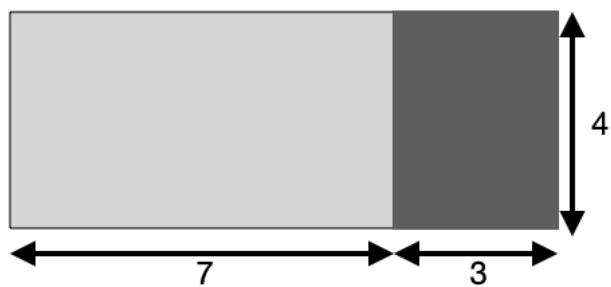

Je vais retourner les rectangles et compter tous les carreaux. Combien vais-je trouver de carreaux ?

Problème 12

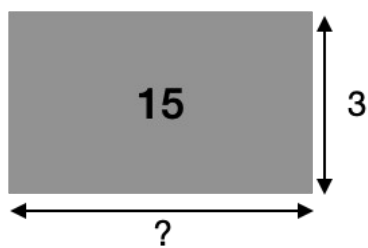

## II. EFFECTUER UNE VALIDATION MATÉRIELLE APRÈS AVOIR RÉSOLU LE PROBLÈME

Pour chacun des douze problèmes proposés, quand la recherche est terminée, le matériel permet de confronter les réponses fournies à la réalité : j'ouvre la boite et je compte les billes, je pose les objets sur la balance, je retourne les bandes et je compte les carreaux, j'ouvre le porte-monnaie et je compte les euros… Dit autrement, ces problèmes permettent une validation par le matériel (Brousseau, 1998). Ainsi, ce n'est ni l'enseignant, ni la majorité qui fournit la réponse correcte, mais l'observation des faits qui fournit la vérité. Pour être satisfaisante, une réponse doit être vraie (et pas seulement bonne, correcte, juste, exacte ou attendue). Qu'est-ce que cela apporte ?

### 1. Convaincre sans juger ni développer de sentiment d'injustice

Un premier bénéfice est d'ordre psychologique. Imaginons un élève qui, pour le problème 3, répond que la bande verte contient 5 cases (peut-être parce qu'il a remarqué qu'elle est longue à peu près comme la moitié de la bande « 10 »). Si on retourne la bande et qu'il constate qu'il s'est trompé, il sera probablement déçu, mais son échec est aussi temporaire et incontestable qu'une chute de vélo : une fois surmonté le premier moment de déception, cela donne surtout envie de faire mieux la prochaine fois.

En revanche, si une correction est effectuée oralement, sans retourner la bande, il n'est pas certain que la réponse fournie et son explication soient plus convaincantes pour l'élève que son propre cheminement. Le risque existe qu'il trouve injuste qu'on félicite ses camarades ayant répondu 6 et pas lui, alors qu'il a fait le même travail : réfléchir et répondre à la question posée par l'enseignant.

Ce sentiment d'injustice, s'il est répété, risque d'entrainer le découragement, le renoncement : de toute façon, je suis nul en maths. L'habitude de chercher à satisfaire l'enseignant en devinant ses attentes, celle de vouloir fournir une réponse « bonne », c'est-à-dire appréciée par l'enseignant, plutôt que vraie, trouvent probablement leur source dans ce type d'expérience. Quand un problème est posé en décrivant oralement ou par écrit une situation fictive, une telle validation n'est généralement pas possible.

## 2. Placer dès le début de la scolarité la question de la vérité au centre de l'activité mathématique

À ce niveau de scolarité, comparer systématiquement les réponses proposées à la réalité du matériel contribue à installer un contrat didactique sain : résoudre un problème, c'est essayer de dire quelque chose de vrai à propos d'objets qu'on ne peut pas ou que l'on ne veut pas compter. Il arrive qu'on se trompe, ce n'est pas grave. On essaie alors de comprendre où on s'est trompé. Parfois, ça aide à réussir les problèmes suivants. Cette question de la vérité est centrale dans la pratique des mathématique quel que soit le niveau (Margolinas, 1993).

Après quelques mois de travail, l'enseignant pourra commencer à poser la question suivante à la fin de chaque problème : « Faut-il vraiment ouvrir la boite (ou retourner la bande, ou vider le porte monnaie…) pour savoir si ce qu'on a trouvé est vrai ? ». Quelques semaines ou mois plus tard, arrivera le moment où, pour certains problèmes, la plupart des élèves répondront que ce n'est pas nécessaire. Cela témoignera d'un apprentissage : dans leur esprit, le raisonnement mathématique utilisé pour répondre à la question est suffisamment consistant et convainquant pour ne pas avoir besoin d'une confirmation matérielle.

L'enseignant pourra alors ajouter progressivement aux problèmes permettant une validation matérielle quelques problèmes portant sur une situation fictive, pour lesquels une validation matérielle est difficile ou impossible.

## 3. Rendre plus efficace le travail sur les procédures

La vérification matérielle permise par les problèmes de notre échantillon, a selon nous un autre effet positif : elle augmente la qualité d'attention des élèves et la clarté des explications de l'enseignant.

D'abord, les élèves sont plus attentifs aux explications de l'enseignant car, sachant déjà si leur réponse est vraie ou non, ils ne sont pas focalisés sur la question « est-ce que j'ai bon ? ». Par ailleurs, les explications de l'enseignant sont plus claires car elles ont une unique fonction : analyser les procédures pour expliquer en quoi certaines sont erronées, d'autres sont efficaces… dans le but d'aider les élèves à résoudre les problèmes à venir. En revanche, après un problème portant sur une situation fictive décrite par un texte, l'enseignant qui cherche à expliciter les procédures doit simultanément indiquer quelle est la réponse correcte, convaincre qu'elle est correcte et que les autres sont erronées. Ces fonctions très différentes sont difficiles à concilier.

# III. TOUJOURS PERMETTRE DE RÉPONDRE « JE NE SAIS PAS ».

Si, dans la rue, je demande comment me rendre à la boulangerie la plus proche, je m'attends à une réponse ressemblant à ceci : « Prenez la deuxième rue à gauche, la boulangerie sera à 100 mètres sur votre droite. ». Je sais que la personne interrogée peut commettre une erreur : oublier l'existence d'une petite rue (la boulangerie est en réalité dans la troisième rue à gauche), confondre droite et gauche ou indiquer un emplacement où il n'y a plus de boulangerie depuis quelques mois). En revanche, je n'envisage pas que la personne interrogée fronce les sourcils et fasse semblant de réfléchir avant de formuler sa réponse alors qu'elle ignore où se trouve la boulangerie. Un contrat implicite dans cette situation est que la personne qui répond cherche à dire la vérité, même si elle n'y parvient pas toujours.

Il nous semble qu'un tel contrat doit être installé dès que possible dans la classe de mathématiques, c'est un premier pas vers ce qu'on appelle en mathématiques la rigueur : quand on émet une affirmation, on a des raisons qu'on juge solides de penser qu'elle est vraie, sinon on ne dit rien, ou on dit qu'on ne sait pas.

Or, en encadrant des groupes d'étudiants préparant l'épreuve de mathématique du concours de recrutement de professeurs des écoles, nous avons pu constater que, même hors d'un contexte d'évaluation, nombre d'entre eux ne voyaient pas d'inconvénient à dire ou écrire une affirmation (résultat intermédiaire ou réponse à la question posée) sans être convaincus que cette affirmation était vraie. Que ces étudiants, par ailleurs intelligents, cultivés, sérieux et travailleurs, considèrent comme acceptable, en mathématiques, d'affirmer quelque chose sans penser que c'est vrai est sans aucun doute un symptôme majeur de leurs difficultés avec cette discipline. C'est en même temps un obstacle qu'il faut lever pour leur permettre de progresser.

Dans les autres disciplines, le critère de réussite peut ne pas être la vérité (arts, EPS) ou bien la vérité peut être plus difficile à établir qu'en mathématiques élémentaire (problème des sources en histoire, démarche expérimentale en sciences). C'est selon nous une raison de plus de s'attacher, en mathématiques, à être exigeant sur la recherche de la vérité.

Nous faisons l'hypothèse qu'en prévoyant de façon systématique le droit de répondre « je ne sais pas » à une question mathématique on peut limiter l'apparition de la forme d'aliénation que nous venons de décrire. C'est par ailleurs un changement qui ne demande aucun moyen particulier : il suffit d'être persuadé que c'est utile pour le mettre en œuvre. Considérons ainsi les problèmes où l'on demande de comparer la longueur de deux grandes bandes (par exemple problème 4). L'enseignant peut demander, après le temps de recherche nécessaire, de répondre en levant un carton parmi les quatre dont chaque élève dispose :

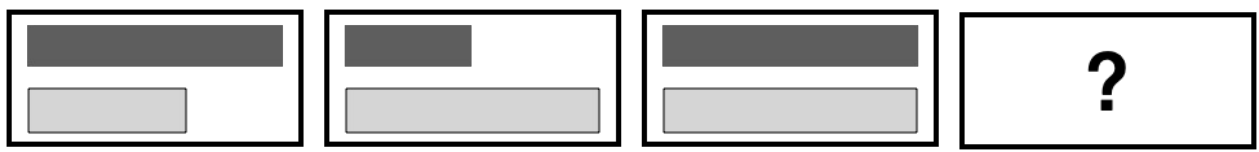

Le quatrième carton évite à l'élève qui ne sait pas de se sentir obligé de répondre un peu au hasard pour obéir à l'injonction et montrer qu'il a travaillé.

Ces réponses « au hasard », qu'on obtient fréquemment si les élèves ne se sentent pas autorisés à répondre qu'ils ne savent pas, sont regrettables par elles-mêmes. Elles ont, en outre, une conséquence tout aussi regrettable : pousser l'enseignant à gaspiller du temps et de l'énergie à tenter d'analyser et de corriger ce qu'il prend à tort pour une erreur. Si un élève répond à tort que la bande rouge est plus longue, non pas parce qu'il le pense, mais parce qu'il faut bien répondre quelque chose, il n'y a rien à analyser. Comme il est souvent difficile de distinguer une véritable une erreur d'une réponse fausse « au hasard », le risque existe aussi de ne pas consacrer suffisamment d'attention à l'analyse de vraies erreurs, ce qui est pourtant essentiel dans le processus d'apprentissage.

Par ailleurs, le fait de pouvoir répondre « je ne sais pas » permet à l'enseignant d'avoir une vue plus réaliste des capacités de ses élèves. Si, pour chaque problème quelques élèves seulement répondent qu'ils ne savent pas, et si ce ne sont pas toujours les mêmes, tout va bien. En revanche si certains élèves répondent presque systématiquement qu'ils ne savent pas, ou si pour un problème donné presque toute la classe ne sait pas, l'enseignant reçoit des indications précieuses pour son action à venir.

## IV. DES PROBLÈMES SANS ÉNONCÉ ÉCRIT, OU PRESQUE, ET UN CONTEXTE RÉCURRENT

On observe fréquemment des séances en principe consacrées à la résolution de problème, mais où l'essentiel du temps est consacré à lire l'énoncé, à expliquer les mots difficiles (voir par exemple Coppé, Houdement, 2002). On y fait probablement un travail utile du point de vue de la langue, mais fort peu de mathématiques. En début de CP, comment permettre la dévolution d'un problème à des élèves qui ne sont pas encore lecteurs ? Cette question est cruciale et elle le reste d'ailleurs au cours de l'année alors que tous les élèves ne progressent pas de la même façon dans l'apprentissage de la lecture. Les problèmes sans énoncé écrit, comme ceux nous avons présentés dans la première partie de l'article, présentent de ce point de vue plusieurs avantages.

### 1. Une passation de la consigne et une appropriation de la question facilitées

L'appui sur un contexte matériel présent dans la classe, sous les yeux des élèves, plutôt que sur un texte — même illustré comme on le voit dans beaucoup de manuels — permet une passation rapide de la consigne.

De plus, l'utilisation récurrente d'un matériel rapidement connu des élèves permet de réduire le temps d'appropriation de la question à presque rien. Les bandes quadrillées utilisées pour les problèmes 3 et 4 en sont un bon exemple. Quand le matériel est connu, les problèmes où on cherche la longueur d'une bande, comme le problème 3, se posent en plaçant les bandes au tableau, sans commentaire ni indication écrite, comme vous l'avez peut-être relevé dans la première partie de l'article. Les problèmes de comparaison de longueur comme le problème 4 peuvent également être ritualisés et posés sans aucun accompagnement de texte : il suffit de placer au tableau des bandes de deux couleurs dont les nombres de cases sont tous connus.

### 2. Un accès aux éléments pertinents du problème facilité

Julo (1995, 2002) indique qu'il n'y a pas deux phases successives de compréhension du problème puis de résolution, que les deux processus sont étroitement liés. Toutefois, il nous semble qu'on peut distinguer un premier niveau de compréhension élémentaire préalable à partir duquel peut s'enclencher le processus par lequel tous les éléments significatifs du problème s'éclairent en même temps que la solution s'élabore. Par exemple, pour le problème 4, surtout si d'autres problèmes de comparaison ont déjà été rencontrés, l'idée qu'on cherche à prévoir quelle bande sera la plus grande quand on effectuera les assemblages, est immédiate. En revanche, le fait qu'il y a une bande rouge et une bande bleue de même longueur est un élément du problème qui peut ne pas apparaître immédiatement et qui fait partie du processus d'élucidation-résolution.

Quand un problème est posé à l'aide d'un texte à des élèves encore peu autonomes en lecture, l'étayage de l'enseignant permet éventuellement à tous les élèves d'atteindre un premier niveau de compréhension : on cherche combien de billes possède Nathalie. En revanche, tous les éléments de l'énoncé ne sont pas intégrés dès la première lecture, même assistée : les billes de Nathalie sont-elles une partie d'un ensemble de billes plus vaste ? Sont-elles obtenues en regroupant des sous-ensembles ?.

Pour intégrer tous les éléments pertinents du problème et donc construire le problème, plusieurs lectures sont généralement nécessaires. Pour les élèves les plus avancés en lecture, l'étayage de l'enseignant en début de séance suffit : ils peuvent ensuite relire le texte pour retrouver une information oubliée ou dont l'importance leur avait échappé. Au contraire, les élèves les plus fragiles en lecture sont perdus dès qu'ils oublient une information (ce qui est presque inévitable) parce qu'ils peuvent difficilement la retrouver dans l'énoncé. De plus, pour ces mêmes élèves, la phase de lecture de l'énoncé même étayée par l'enseignant demande un effort très important. Il est possible qu'ils n'aient plus la capacité à se concentrer encore pour répondre à la question posée.Ainsi, les problèmes à texte utilisés prématurément peuvent mettre artificiellement en échec en mathématiques des élèves dont les seules difficultés réelles concernent la lecture.

## 3. Une présentation unique mais des problèmes variés

La possibilité d'aborder des problèmes plus nombreux en supprimant la phase de lecture du texte n'aurait pas d'intérêt si la répétition de problèmes dans une présentation identique conduisait à des tâches routinières, sans réel contenu mathématique. Il n'en est rien avec les problèmes que nous proposons : l'uniformité de la présentation n'empêche pas la variété des situations mathématiques, des niveaux de difficulté et des procédures envisageables. En fait, à partir d'un matériel récurrent nous installons un contexte mathématique source d'une variété de problèmes. Comme nous l'illustrons ci-dessous.

*Problèmes 13 à 18 (recherche de la longueur d'une bande)*

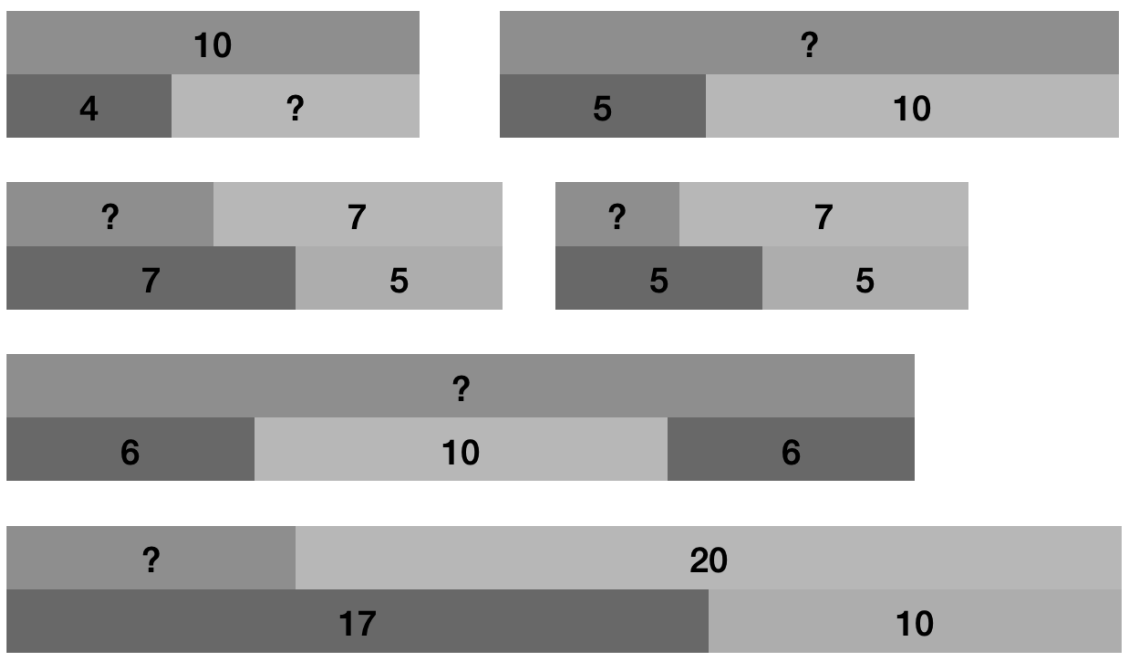

*Problèmes 19 à 21*

Si l'enseignant précise, quand il y a plusieurs bandes inconnues, que ces bandes sont identiques,

la variété des problèmes s'accroit encore sans que la difficulté de compréhension de la question augmente.

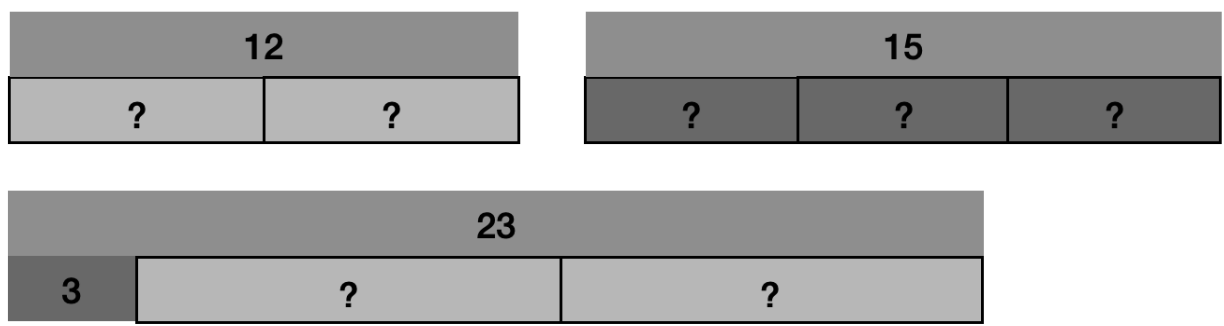

*Problèmes 22 et 23 (comparaison de la longueur de deux bandes obtenues par assemblage)*

Par ailleurs ce matériel permet de poser des problèmes de comparaison qui, à notre connaissance, sont très peu proposés dans les écoles : je fais deux grandes bandes en mettant bout-à-bout les bandes de même couleur. Laquelle sera la plus longue ?

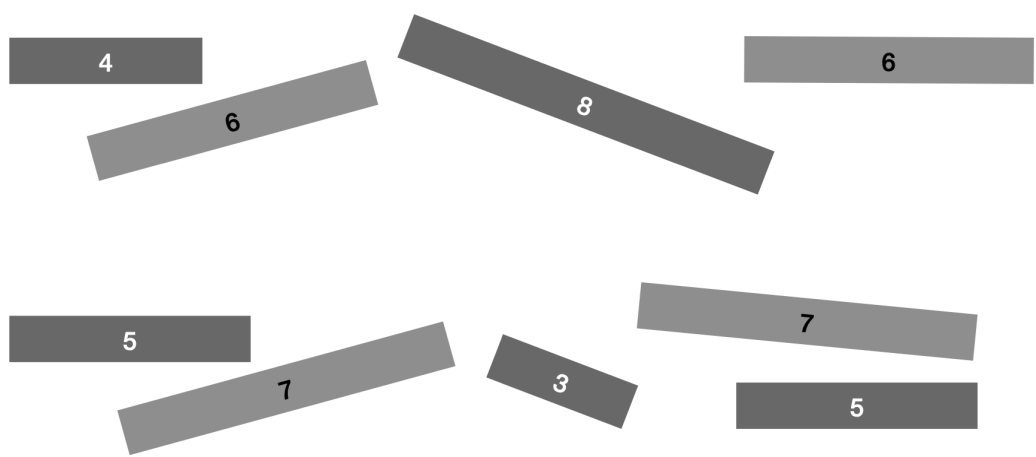

On peut ainsi proposer des problèmes sans texte (ou avec très peu de texte avec d'autres supports que les bandes) qui permettent à toute la classe de se confronter à un grand nombre de situations mathématiques variées, puisque le temps consacré à la lecture de l'énoncé n'est plus nécessaire.

## 4. Une dévolution facilitée

Nous avons détaillé les problèmes portant sur des bandes qui sont pour nous emblématiques : ritualisés, ils permettent de se passer entièrement de texte écrit ou oral. Il en va de même pour les problèmes portant sur les cartes à points ou les rectangles quadrillés. Les autres problèmes proposés au début de cet article ne sont pas tous aussi radicaux de ce point de vue. Certains nécessitent un accompagnement oral que nous ne détaillerons pas ici, mais qui peut être rendu assez léger si l'enseignant choisit de poser un certain nombre de problèmes à propos du même type de matériel.

Remarquons que, dans les situations de la vie courante où l'on rencontre un problème mathématique, celui-ci n'est pas posé par un texte. La question est généralement ce qui apparait en premier : ai-je assez d'argent pour m'acheter ce livre et aller au cinéma ? Ce meuble peut-il passer dans l'escalier sans que je le démonte ? À quelle heure faut-il que je parte pour arriver chez Grand-mère avant 18 h ? On collecte ensuite les informations nécessaires pour répondre à la question. Dans les problèmes sans texte, comme dans ces situations de la vie courante, la question est posée d'emblée, ce qui facilite la dévolution.

## 5. Un tremplin pour la résolution de problèmes avec énoncés textuels

Il faudra bien que les élèves abordent un jour des problèmes posés par un énoncé écrit puisque tel est l'usage. À cet égard, proposer des énoncés évoquant les situations travaillées auparavant avec du matériel nous semble une première étape possible. On pourra ainsi proposer par écrit le problème suivant :

*Sur une balance, je place :*

*15 billes rouges et 12 billes bleues à gauche,*

*9 billes jaunes et 20 billes vertes à droite.*

*Quel côté va descendre ?*

Nous pensons que la fréquentation de nombreux problèmes tels que nous les avons décrits plus haut est plus une aide qu'un obstacle pour aborder ces problèmes à énoncé écrit. Nous manquons de recul, d'observations conduites avec rigueur, mais le retour d'un enseignant ayant utilisé nos propositions pendant une année de CP apparaît tout de même encourageant : conscient d'être en décalage avec les méthodes utilisées dans les autres CP de l'école, il n'a pas voulu participer à l'élaboration d'une évaluation de fin d'année commune à tous les CP. Ses élèves ont donc été confrontés à des problèmes à énoncé alors qu'ils n'en avaient rencontrés que quelques uns en toute fin d'année. Ils ont pourtant eu les meilleurs résultats, tant en calcul qu'en résolution de problèmes.

# V. DÈS LE DÉBUT, CALCULER ET NON COMPTER

Comme Brissiaud (2003) l'a souligné, s'appuyer essentiellement sur le comptage n'est pas une voie d'accès vers le calcul, mais un obstacle à l'apprentissage du calcul. C'est pourtant ce qui est proposé dans le guide « pour enseigner les nombres, le calcul et la résolution de problèmes au CP » (MEN, 2021). Par ailleurs, calculer ne signifie pas forcément choisir une opération pour l'effectuer. Nous donnons dans cette partie quelques procédures de calcul réfléchi pouvant être utilisées au CP avant l'introduction des opérations.

## 1. Une proposition qui interpelle

Le guide pour le CP (MEN, 2021) distingue trois stratégies de résolution de problèmes additifs : 1/stratégies de dénombrement plutôt élémentaires ; 2/stratégies de dénombrement s'appuyant sur des représentations symboliques des collections ; 3/stratégies de (ou proches du) calcul, plus ou moins explicitées ou formalisées. L'objectif désigné pour la résolution de problème au cycle 2 (MEN, p. 87) est « d'arriver progressivement à des procédures relevant de la stratégie 3 et en particulier à la production d'écritures mathématiques. ». Pour cela, les enseignants sont invités à repérer les stratégies utilisées par leurs élèves et à passer progressivement de la représentation[6] à la modélisation du problème utilisant des schémas en barre. La modélisation est présentée comme une aide à la résolution de problème. La progression proposée consiste à « passer du dessin figuratif au schéma grâce au matériel » (ib., p. 92), en association avec un jeu sur les valeurs numériques (ib., p. 94) :

> *Il est important de faire travailler les élèves d'abord dans le cadre de problèmes additifs simples (de type partie 1 + partie 2 = tout) sur des nombres de petites tailles, en lien avec l'apprentissage de la numération et de l'addition, pour construire le sens de l'opération et installer des automatismes ainsi que le lien naturel entre les nombres.*
>
> *En s'appuyant sur le matériel manipulé en maternelle (cubes emboîtables), les problèmes de parties-tout se modélisent progressivement avec des schémas en barres. Les cubes (de couleur) emboîtés deviendront, par un travail d'appropriation pas à pas, les barres rectangulaires dans un schéma.*

---
[6] Dans le guide la différence entre « représenter » et « modéliser » est explicitée de la façon suivante (MEN, p. 89) :

L'extrait qui figure en annexe 1 (ib., p. 94) illustre une résolution du problème « Léo a 7 billes rouges et 5 billes bleues. Combien Léo a-t-il de billes en tout ? ». Il nous parait significatif de la démarche d'abstraction progressive développée dans le guide.

Résoudre le problème des billes de Léo comporte, pour un élève de CP, deux difficultés. Il faut d'abord interpréter l'énoncé, le reformuler éventuellement, intégrer que cela revient à se demander « 7 trucs et encore 5 trucs, c'est combien de trucs ?». Il faut ensuite répondre à cette question. Le document qui figure en annexe 1 n'apporte pas d'aide pour la première difficulté : répondre à la question en prenant 7 cubes rouges et 5 cubes bleus puis en les dénombrant suppose qu'on a interprété correctement l'énoncé. Il ne décrit pas non plus la façon dont les élèves répondent à la question « 7 et encore 5, c'est combien ? ». À l'exception de la dernière étape, les élèves peuvent toujours dénombrer les objets ou les cases du dessin. On maintient ainsi longtemps les élèves dans une procédure qui est précisément celle qu'on souhaite qu'ils abandonnent. Par ailleurs, on voit bien comment, du point de vue graphique, chaque étape rapproche un peu du schéma en barre, mais cette évolution dans la nature du support matériel ou visuel proposé aux élèves ne nous semble pas être une abstraction significative. Du point de vue mathématique, nous voyons plutôt une continuité, puis une rupture à la dernière étape : pour la première fois, on ne peut plus compter.

Le guide (MEN, 2021) propose aussi de faire le lien avec la numération décimale pour obtenir le résultat souhaité, mais, les représentations alors proposées renvoient encore au comptage (voir pages 89 – 92). L'usage de cette seule stratégie risque selon nous de renforcer ce que l'on veut éviter.

## 2. Quelle(s) alternative(s) ?

Pour éviter de renforcer chez les élèves le recours systématique au dénombrement, nous proposons d'énoncer aux élèves la règle suivante :

*Quand on résout un problème, on ne compte qu'à la fin, quand on a répondu. On compte pour savoir si on a trouvé la vérité ou si on s'est trompé.*

L'enseignant qui accepte cette proposition n'encouragera pas ses élèves à prendre des cubes pour

---

« Représenter, c'est traduire par un dessin ou un schéma la situation. Le fait de représenter la situation permet de l'appréhender et de favoriser l'entrée dans la résolution. Certaines représentations (souvent de type pictural) ne sont pas traduisibles par un calcul. » ; « Modéliser, c'est traduire mathématiquement la situation. La modélisation amène ensuite à la procédure et au calcul ; elle rend la réalité calculable. Il s'agit d'un processus qui peut prendre appui sur diverses représentations. » résoudre le problème de Léo. Au contraire, il dira à celui qui prend l'initiative de le faire « Tu es vraiment sûr que tu veux compter les cubes ? Pour réussir, il faut trouver combien Léo a de cubes sans les compter ? »

Imaginons que, suite à un entraînement préalable sur des problèmes analogues, un élève représente par un schéma en barre le problème des billes de Léo avant de le résoudre, comme ceci :

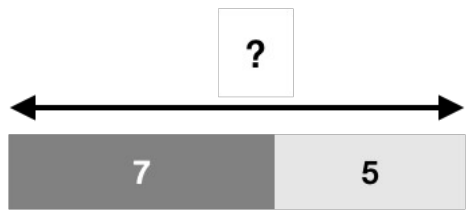

Ce schéma témoigne que cet élève a bien compris la structure du problème, mais il ne l'aide pas à répondre à la question « 7 billes et encore 5 billes, c'est combien de billes ? ».

Cette difficulté n'est pas propre au schéma en barres. Imaginons que, dans une autre classe de CP, on ait appris tôt dans l'année que, quand on réunit une collection de 7 éléments et une de 5 éléments, le nombre d'éléments en tout se note 7+ 5. Un élève de cette classe écrira peut-être que le nombre de billes de Léo est 7 + 5, mais cela ne l'aidera pas à déterminer que 7 + 5, c'est la même chose que 12. L'écriture 7+5 = 12 est seulement, comme le schéma en barre à la fin de la l'annexe 1, une façon de traduire *a posteriori* ce qu'on a fait. La question qui se pose à l'élève qui vient de réaliser le schéma en barre ou d'écrire que le nombre cherché est 7 + 5 est exactement celle qu'il a pu se poser avant de réaliser le schéma ou l'écriture symbolique : 7 et encore 5, c'est combien ?

Comment un élève peut-il déterminer que 7 et encore 5, c'est 12 en début de CP ? Répondre à cette question dépasse le cadre que nous nous sommes fixé ici : décrire quelques conditions favorisant l'apparition d'un contrat sain pour la résolution de problèmes numériques. Néanmoins, il nous semble difficile d'encourager les enseignants à dissuader leurs élèves de dénombrer sans esquisser une alternative. Il se peut que, pour quelques élèves, le fait numérique « 7 et encore 5, c'est 12 » soit déjà mémorisé et facile à mobiliser en début de CP. Ces élèves constituent une minorité pour laquelle nous ne nous inquiétons pas. Pour les autres, si l'enseignant ne donne pas explicitement d'exemples de procédures de calcul, le recours au dénombrement est la seule option disponible. Donnons un exemple de procédure de calcul envisageable à ce stade en s'appuyant sur les « cartes à points » recto-verso (figure 1). Chaque carte comporte sur sa face cachée des points disposés en utilisant les constellations du dé (le nombre de points étant évidemment celui écrit au recto).

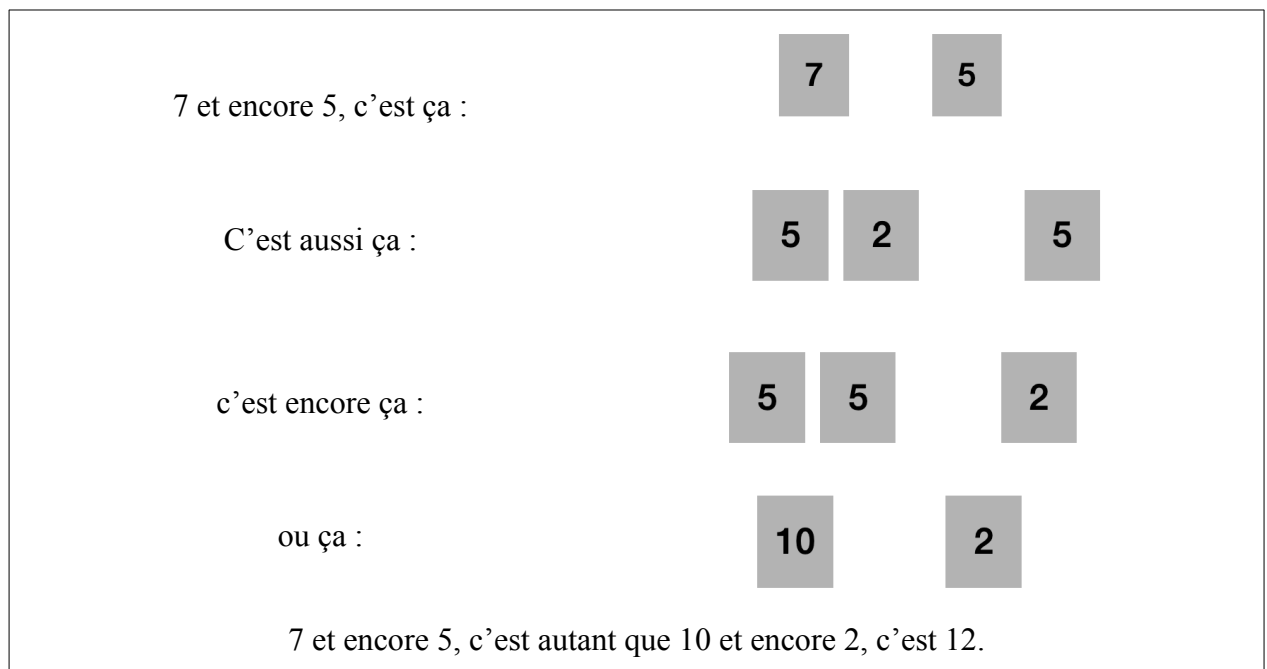

**Figure 1** : une première décomposition possible

Cette procédure de calcul mobilise les faits numériques suivants :

    7, c'est 5 et encore 2.

    5 et encore 5, c'est 10

et ces propriétés plus générales :

    Si on déplace des points, leur nombre ne change pas

    Si on partage un groupe de points en deux groupes plus petits, ou si on rassemble deux petits groupes pour faire un grand, le nombre de points ne change pas.

    10 et encore un petit nombre, ça s'écrit avec 1 suivi du petit nombre.

Nous présentons en figure 2 un autre exemple de procédure pour le même calcul. Dans cette seconde procédure, les faits « 7, c'est 6 et encore 1 » et « 6, c'est 5 et encore 1 » peuvent être évoqués séparément ou vus comme des cas particuliers de la propriété « le nombre qui suit dans la comptine, c'est un de plus ».

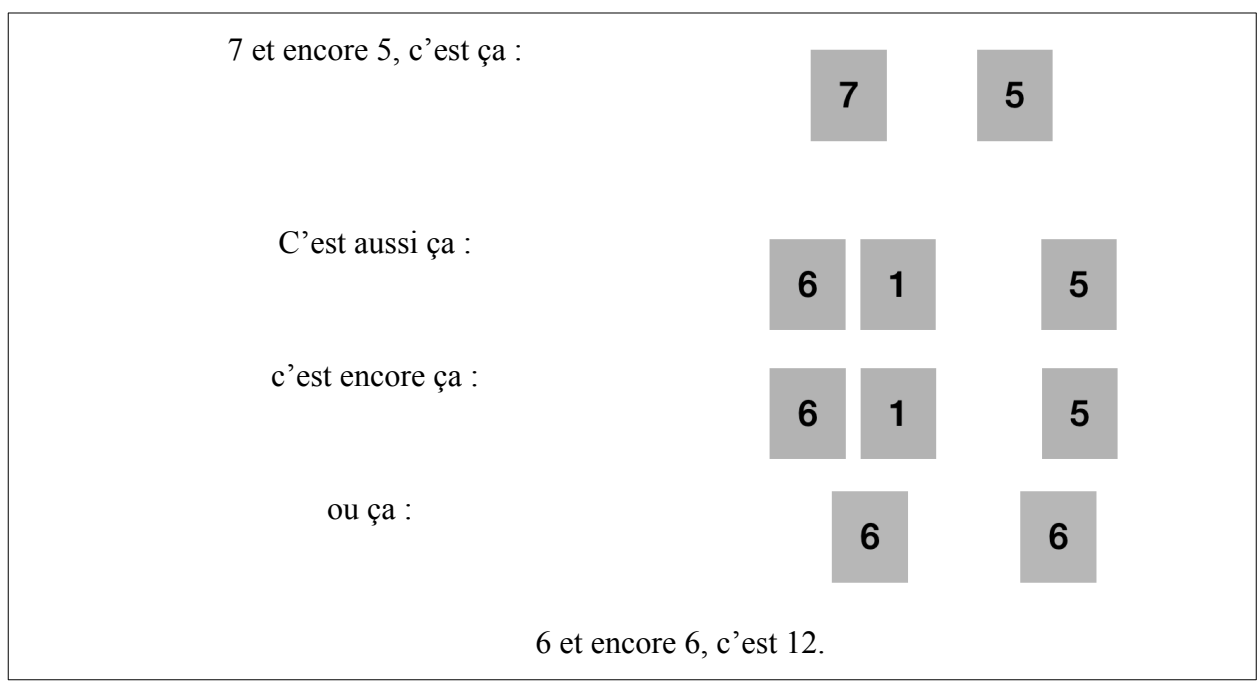

**Figure 2 :** une autre décomposition possible

Ces procédures ne sont pas spontanées, les faits numériques et les propriétés mobilisés résultent d'un apprentissage structuré des premiers nombres. Comment sont-ils enseignés et appris ? Il y a bien sûr plusieurs réponses envisageables, mais quelle que soit la réponse qu'on donne à cette question, construire un répertoire de faits numériques et de propriétés des nombres afin que les élèves puissent, comme dans les exemples ci-dessus, agir sur les nombres autrement qu'en comptant nous semble une priorité. C'est une façon de permettre ce que Butlen (2007) appelle l'alternance entre sens et technique.

## CONCLUSION

Nous avons développé quatre conditions didactiques qui visent à installer un contrat didactique pertinent pour la résolution de problèmes numériques, dès le début du cycle 2, en amont d'un travail sur la représentation des problèmes : effectuer une validation matérielle ; permettre de répondre « je ne sais pas » ; proposer des problèmes sans énoncé écrit ou presque, avec un contexte récurrent ; demander de calculer sans compter/ dénombrer.

La mise en place de ces conditions en classe est plus ou moins aisée, et plus ou moins liée à un choix de problèmes. Ainsi, permettre systématiquement (et le rappeler régulièrement) de répondre qu'on ne sait pas ne demande aucun moyen particulier... il suffit que l'enseignant soit convaincu que c'est utile pour qu'il le fasse. En revanche, poser des problèmes portant sur un matériel présent dans la classe et en utilisant peu ou pas de texte écrit suppose un travail plus important. De tels problèmes sont peu présents dans les diverses sources disponibles et il faudra peut-être les construire. Par ailleurs, un travail de préparation matérielle est souvent nécessaire : se procurer une balance de Roberval, fabriquer les bandes quadrillées à afficher au tableau... Il ne faudrait cependant pas exagérer ces difficultés : les bandes quadrillées qu'on aura fabriquées serviront toute l'année et la balance de Roberval dénichée avec difficulté pourra servir toute une carrière. Enfin, expliciter que résoudre un problème c'est répondre sans compter est sans doute celle de nos quatre propositions qui exige le changement le plus radical des pratiques, car de nombreuses ressources s'appuient au contraire fortement sur le dénombrement, en cohérence avec les préconisations du guide publié par le MEN en 2021.

## RÉFÉRENCES BIBLIOGRAPHIQUES

# ANNEXE 1
# EXTRAIT DU GUIDE

UN EXEMPLE DE PROBLÈME ET DE MODÉLISATION PROGRESSIVE

PAR LE SCHÉMA EN BARRES

➔ « *Léo a 7 billes rouges et 5 billes bleues. Combien Léo a-t-il de billes en tout ?* »

La résolution de ce problème à l'aide de 7 cubes rouges :

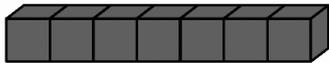

et 5 cubes bleus :

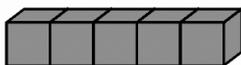

fait apparaître l'assemblage :

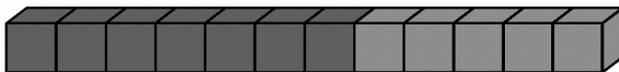

puis le schéma :

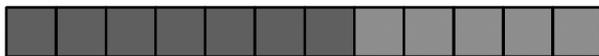

et enfin le schéma en barres :

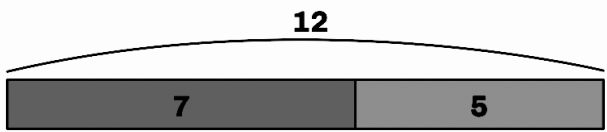

*figure 1*

----